\documentclass{birkjour}
 \newtheorem{thm}{Theorem}[section]
 \newtheorem{cor}[thm]{Corollary}

 \newtheorem{rem}[thm]{Remark}
 
 \numberwithin{equation}{section}

\newcommand{\R}{{\mathbb R}}

\newcommand{\D}{{\mathbb D}}

\renewcommand{\k}{{\bf k}}

\newcommand{\e}{{\bf e}}
\newcommand{\edag}{{\bf e^\dagger}}

\newtheorem{Pa}{Paper}[section]
\newtheorem{Tm}[Pa]{{\bf Theorem}}

\newtheorem{Cy}[Pa]{{\bf Corollary}}

\newtheorem{Dn}[Pa]{{\bf Definition}}

\begin{document}

\title[Kolmogorov's  axioms...]{Kolmogorov's axioms  for probabilities with va\-lues in  hyperbolic numbers}

\author[D. Alpay]{Daniel Alpay}
\address{(DA) Department of mathematics,
Ben-Gurion University of the Negev, P.O. Box 653, Beer-Sheva
84105, Israel}

\email{dany@math.bgu.ac.il}

\author[M. E. Luna--Elizarrar\'as]
{M. Elena Luna--Elizarrar\'as}

\address{(ME L-E) Escuela Superior
de Fisica y Matem\'aticas, Instituto Polit\'ecnico Nacional,
Mexico City, M\'exico.}

\email{eluna@esfm.ipn.mx}

\author[M. Shapiro]{Michael Shapiro}

\address{(MS) Escuela Superior
de Fisica y Matem\'aticas, Instituto Polit\'ecnico Nacional,
Mexico City, M\'exico.}

\email{shapiro@esfm.ipn.mx}

\date{}

\begin{abstract}
We introduce the notion of a probabilistic measure which takes values in hyperbolic  numbers and which satisfies the system  of axioms  generalizing directly   Kolmogorov's  system of  axioms.  We show  that this new measure verifies  the usual  properties   of  a  probability;  in  particular,  we treat  the  conditional  hyperbolic  probability  and we prove the hyperbolic analogues of the multiplication  theorem, of the law of total  probability and of  Bayes'  theorem.  Our probability  may take values  which are zero--divisors and  we discuss carefully this peculiarity.
\end{abstract}

\subjclass{} \keywords{} \maketitle


\date{}

\maketitle

\section{Introduction}

The hyperbolic numbers (called also split--complex, or perplex, or double numbers, etc.) are known since long ago but they are not as popular as complex numbers or quaternions. At the same time they possess many  interesting properties; in particular, the  ring $\D$ of hyperbolic numbers admits a partial order  $\preceq$ which has a good compatibility with the other algebraic  structures of $\D$.  Consider the inequality  $ 0 \preceq x \preceq 1$; it turns out that it has a well--defined  set of solutions in $\D$ and one can think of them as of the probabilities of some random events.

This was a motivation of the present work: to test how this conjecture operates. First of all, we give  a review of  hyperbolic numbers making a special  emphasize  on the properties of  non--negative  hyperbolic numbers. Next,  we introduce  direct  generalizations of  Kolmogorov's  axioms where  a  probabilistic  measure takes  values in hyperbolic numbers. It is followed  by a series  of the immediate  properties of such probabilistic  measures.  The last Section 5  ``Conditional probability"  introduces this notion, including the case of probabilities  which  are  zero--divisors  in  $\D$,  and  presents  the  hyperbolic  generalizations of several  classic  facts:  multiplication theorem, independence of random events, law of total probability,    Bayes'  theorem.

Altogether, we  have shown  that the basic  facts  (elementary but underlying)  of the  classic  probability  theory extend onto the situation under consideration.  Thus one can expect that the whole  building can be constructed  as well.

Our  approach  can be seen in different ways.  First of all, since we replace the range $\R$ of probabilistic measures  by a hypercomplex system $\D$  then the approach can be interpreted as an attempt to consider probabilistic measures with non--numerical  (in the sense of non--real)  values.  The most  renowned research line in this direction is that of quantum probability but the two  are rather  distant  from each other, see, e.g.,  \cite{Accardi}, \cite{Parth}, \cite{Meyer}. 

On the other hand, the hyperbolic numbers can be seen as a real two--dimensional  algebra with the underlying linear  space  $\R^2$. Hence, the hyperbolic probability  has the following interpretation: one deals with a stochastic experiment  which generates  the necessity to endow  the $\sigma$--algebra of the events  with two probabilistic measures  which are seen  as $\R^2$--valued measures;  what is more,  a rich multiplicative structure  is introduced on the  range of such measures. 

Such situations may emerge in mathematical statistics in testing composite hypotheses.

Another example is provided  by thermodynamics and statistical physics. Consider a physical system which has  two (or more)  minima of free energy. If  the system is in an equilibrium  then it can be in any of these states with certain probabilities but it cannot be known  for sure  in which of them;  this is exactly the situation we are interested in.  

We believe that our approach will be useful in  treating such situations although in the present work we limit ourselves with considering the basics of purely mathematical theory.

The Mexican authors were partially supported  by Instituto Polit\'ecnico Nacional in the framework of COFAA and SIP programs,  as  well  as  by  SNI--CONACYT. \\


\section{A review of  hyperbolic numbers}
Information about  hyperbolic numbers is dispersed in many sources. We concentrate in this section some basic facts which can be found in more details in  \cite{uno}, \cite{cinco}.  \\

The ring of hyperbolic numbers is the commutative ring 
 $\D$  defined as    
$$\D := \left\{ a+b \k \mid a,\, b\in\R ;\; \k ^{2} =1,\, \k \notin \R \right\}.$$

There is a  conjugation, the $\dagger$--conjugation, on hyperbolic  numbers:
$$ \mathfrak z^{\dagger} := a - b \k .  $$

This $\dagger$--conjugation is  an additive, involutive and multiplicative operation on $\D$:
\begin{enumerate}
\item $ (\mathfrak z+ \mathfrak w)^\dagger = \mathfrak z^\dagger + \mathfrak w^\dagger$; 
\item $ (\mathfrak z^\dagger )^\dagger = \mathfrak z $;  
\item $( \mathfrak z \mathfrak w)^\dagger = \mathfrak z^\dagger \mathfrak w^\dagger$.  \\
\end{enumerate}

Note that given  $\mathfrak z= a + b \k \in \D$,  then 
$$ \mathfrak z \mathfrak z^\dagger = a^2 - b^2 \in \R,$$
from which  it follows that any hyperbolic number $\mathfrak z$ with $ \mathfrak z \mathfrak z^\dagger \not= 0$ is invertible, and its inverse is given by 
$$ \displaystyle  \mathfrak z^{-1} = \frac{ \mathfrak z^\dagger}{ \mathfrak z \mathfrak z^\dagger} . $$
If, on the other hand, $\mathfrak z \not= 0$ but $ \mathfrak z\mathfrak z^\dagger =a^2 - b^2 =0$ then $\mathfrak z$ is a zero--divisor. In fact there are no other zero--divisors. We denote the set of zero--divisors by $\mathfrak S_\D$, thus
$$\mathfrak{S}_\D := \left\{\mathfrak z =  a + b \k   \mid  \mathfrak z \neq 0,\, \mathfrak z\mathfrak z^\dagger= a^2 - b^2 =0  \right\}.$$
It turns out that there  are two very special zero-divisors in $\D$. Set 
$$\e :=\frac{1}{2}(1+\k),   $$ 
then its $\dagger$-conjugate is 
$$\edag :=\frac{1}{2} (1-\k ).$$
It is immediate to check that $\e$ and $\edag$ are zero-divisors, and they are mutually complementary idempotent elements.  Thus, the two sets 
$$\D _\e := \e\cdot \D \,\mbox{ and }  \, \D _\edag := \edag \cdot \D$$ 
are (principal) ideals in the ring $\D$ and they have the   properties: 
$$\D _\e \cap \D_\edag =\left\{0\right\}$$ 
and 
\begin{equation}\label{unouno}
\D = \D _\e + \D _\edag .
\end{equation}

Formula (\ref{unouno})  is called the idempotent decomposition of $\D$. Every hyperbolic  number $\mathfrak z= a + b \k $ can be written as 
\begin{equation}\label{unodos}
\mathfrak z =a + b \k = (a+b) \e + (a-b) \edag =:\nu_1 \e + \nu_2  \edag .
\end{equation}
Formula  (\ref{unodos})  is   called the idempotent representation of a hyperbolic number. It  has a remarkable feature: the algebraic operations of addition, multiplication, taking of inverse, etc. can be realized component-wise.  \\

Observe that the sets  $\D_\e$ and $\D_\edag$  can be written as
$$ \displaystyle  \D_\e = \left\{ r \e \mid r \in \R  \right\} = \R \e;    \qquad  \D_\edag = \left\{ t \edag \mid t \in \R  \right\} = \R \edag. $$ 

\begin{rem}\label{useful property}
It will be useful to have in mind the following properties:
\begin{enumerate}
\item[(a)] $\mathfrak z  \in \D_\e$ if and only if $ \mathfrak z \e = \mathfrak z  $;
\item[(b)] $\mathfrak z  \in \D_\edag$ if and only if $ \mathfrak z \edag = \mathfrak z  $.  \\
\end{enumerate}
\end{rem}

The set of non-negative hyperbolic numbers is 
$$\D ^+:= \left\{\nu_1 \e + \nu_2 \edag \mid \nu_1 , \, \nu_2 \geq 0\right\}.$$
We will need two more sets:
$$  \D_\e^+ :=  \D_\e \cap \D^+  \setminus \{0\} \qquad  {\rm and }  \qquad  \D_\edag^+ :=  \D_\edag \cap \D^+  \setminus \{0\} .$$

\medskip

Let us now define on $\D$ the next relation: given $\mathfrak z_1, \mathfrak z_2 \in \D$ we write $\mathfrak z_1\preceq \mathfrak z_2$ whenever $\mathfrak z_2 - \mathfrak z_1 \in \D ^+$; this relation is reflexive, transitive and antisymmetric and therefore it defines a partial order on $\D$. Also, if we take $\alpha,\beta \in \R$ then $\alpha\preceq \beta$ if and only if $\alpha\leq\beta$, thus $\preceq$ is an extension of the total order $\leq$ on $\R$.  \\

The next properties of the order $\preceq$ will be useful in subsequent computations (for more details see \cite{uno}). Let $\mathfrak x, \mathfrak y, \mathfrak z, \mathfrak w \in \D.$
\begin{enumerate}
\item If $\mathfrak x\preceq \mathfrak y$ and $\mathfrak z\in \D ^+$, then $\mathfrak z \mathfrak x\preceq \mathfrak z \mathfrak y.$
\item If $\mathfrak x\preceq \mathfrak y$ and $\mathfrak z\preceq \mathfrak w$, then $\mathfrak x+\mathfrak z\preceq \mathfrak y+ \mathfrak w$.
\item If $\mathfrak x\preceq \mathfrak y$, then $-\mathfrak y\preceq -\mathfrak x$.
\end{enumerate} 

Thanks to the good properties of the partial order  $\preceq$, one defines the hyperbolic--valued  modulus   on  $\D$ by
\begin{equation}\label{def_hyp_norm}
| \mathfrak z |_\k = | \nu_1 \e + \nu_2 \edag |_\k := | \nu_1 | \e + | \nu_2 | \edag \in \D^+ , 
\end{equation}
where $| \nu_1|$, $|\nu_2|$  denote the usual  modulus  of  real numbers.  The subindex  $\k$  is used to emphasize  that this modulus is linked to the hyperbolic  numbers with the imaginary unit $\k$.  Moreover, the name ``hyperbolic--valued modulus"  for  (\ref{def_hyp_norm})   is justified by the following properties (see \cite{uno}, \cite{cinco}):

\begin{enumerate}

\item[({\sc i})]  $| \mathfrak z |_\k =0$ if and only if $\mathfrak z =0$.

\item[({\sc ii})]  $ | \mathfrak w  \mathfrak z |_\k = | \mathfrak w |_\k \cdot | \mathfrak z |_\k $.

\item[({\sc iii})]  $| \mathfrak w +  \mathfrak z |_\k \preceq | \mathfrak w |_\k + | \mathfrak z |_\k$ for any $\mathfrak z$, $\mathfrak w \in \D$.

\end{enumerate}

In particular, one may talk about the supremum  $\sup_\D$  of bounded sets in $\D$ with respect to  this hyperbolic--valued modulus. Indeed, let $A\subset \D$, if there exists $M\in \D ^+$ such that $\left|x\right|_{\k}\preceq M$ for any $x\in A$, we say that $A$ is a $\D$-bounded set. Introduce
\begin{eqnarray*}
A_1 &:=& \left\{ x\in \R \mid \exists  \;  y \in \R,\, x\e +y\edag \in A\right\},\\
A_2 &:=& \left\{ y\in \R \mid \exists \;   x \in \R,\, x\e +y\edag \in A\right\};
\end{eqnarray*}
if $A$ is a $\D$-bounded set then $A_1$ and $A_2$ are bounded, and the $\sup _\D A$ can be computed as 
$${\rm sup} _\D A = \sup A_1 \e + \sup A_2 \edag .$$

It is worth noting that some  hyperbolic  modules can be endowed with a hyperbolic--valued
norm. These norms have the expected properties, that is, if a hyperbolic  module  $W$ has a 
hyperbolic--valued   norm  $\| \cdot \|_\D $, then the  latte satisfies:

\begin{enumerate}

\item  $\| x \|_\D \succeq 0$  for all $x \in W$ and  $\| x \|_\D = 0$ if and only if $x =0\in W$.

\item  $\| \mathfrak z  x \|_\D =  | \mathfrak z |_\k  \| x \|_\D $  for  all  $\mathfrak z \in \D$ and for all $x \in W$.

\item  $ \| x + w \|_\D  \preceq \| x \|_\D  +  \| w \|_\D $  for  all  $x, \, w \in W$.

\end{enumerate}

The strength of hyperbolic--valued norms defined  on  hyperbolic modules has been  exploited in \cite{uno} and  \cite{ocho}. In the latter   a version  of  Hahn--Banach Theorem  for  hyperbolic modules  has  been  proved.  \\


\section{$\D$--valued  probability}

\begin{Dn}
Let  $( \Omega ,\Sigma)$  be  a  measurable  space,  a  function
$$ P_\D:A \in \Sigma \mapsto P_\D(A) \in \D $$
with the properties:
\begin{enumerate}
\item[({\sc i})] $P_\D (A) \succeq 0$ \; $ \forall \, A \in \Sigma$; \\

\item[({\sc ii})] $P_\D (\Omega)= \mathfrak p$, where   $\mathfrak p $ takes one of the three possible values $1$, $\e$,  $\edag$; \\

\item[({\sc iii})]  given a  sequence  $\displaystyle  \left\{ A_n \right\} \subset \Sigma $  of   pairwise  disjoint  events,  then  
$$  \displaystyle  P_\D \left(  \bigcup_{n=1}^\infty A_n \right) =  \sum_{n=1}^\infty P_\D(A_n), $$
\end{enumerate}
is called  a  $\D$--valued  probabilistic  measure,  or  a  $\D$--valued  probability,  on the  $\sigma$--algebra  of  events  $\Sigma$.  The  triplet  $(\Omega, \Sigma ,P_\D)$  is  called  a  $\D$--probabilistic  space.  \\
\end{Dn}

Every  $\D$--valued  probabilistic  measure  can  be  written  as
\begin{equation}\label{cartesian and idempotent representation}
P_\D (A) = p_1 (A)    +  p_2 (A) \k = P_1 (A) \e + P_2 (A) \edag 
\end{equation}
with   $ P_1 (A)  = p_1(A)  +p_2(A); \;  P_2(A) =p_1(A) -p_2(A) $.  
The   property ({\sc i})  of   $P_\D $  implies  that
$$ P_1(A) \geq 0, \quad  P_2(A) \geq 0 \quad  \forall A \in \Sigma.  $$
The  property  ({\sc ii}) gives:
$$  P_\D  (\Omega )=\mathfrak p =  P_1(\Omega) \e + P_2 (\Omega)\edag ,$$
that is:

\begin{enumerate}
\item[(1)] If  $\mathfrak p =1$ then    $ P_1(\Omega) =1$, $P_2(\Omega) =1$. \\

\item[(2)]  If  $\mathfrak p = \e $  then    $ P_1(\Omega) =1$, $P_2(\Omega) =0$. \\

\item[(3)]  If  $\mathfrak p = \edag $  then    $ P_1(\Omega) =0$, $P_2(\Omega) =1$. \\

\end{enumerate}

The  property  ({\sc iii})  leads to 
$$  \displaystyle  P_i  \left( \bigcup_{n=1}^\infty A_n \right) = \sum_{n=1}^\infty P_i (A_n) \quad {\rm for} \quad i=1 \;\; {\rm and} \;\; 2 .$$
Hence, to define a $\D$--valued  probabilistic measure is equivalent to consider, on the same measurable space, a  pair of unrelated, in general, usual $\R$--valued   measures. In the case (1) both $P_1$ and $P_2$ are probabilistic  measures; in case (2)  $P_1$ is a probabilistic measure and $P_2$ is a trivial one; in case (3) $P_1$ is a trivial measure and $P_2$ is a probabilistic measure.  The cases  (2) and (3) can be  seen as two   options of embedding the classic real--valued  probabilistic measures into our new concept of $\D$--valued  probabilistic measures: we identify such real--valued measures  with $\D$--probabilistic  measures  which takes  as  its values  only  zero--divisors.  \\

\section{Properties of $\D$--valued  probabilistic  measures}

\bigskip

\begin{enumerate}

\item[{\bf (I)}] Given $A \in \Sigma$,  then  $P_\D(A) +P_\D (A^C) =  \mathfrak p $  where $A^C \in \Sigma $ is the complement of $A$.

\begin{proof}   $A \cup A^C = \Omega$,  $A \cap A^C = \emptyset$, hence  
$$P_\D (A) + P_\D (A^C) = P_\D (\Omega)= \mathfrak p  .  $$  
\end{proof}

\item[{\bf (II)}]  $P_\D (\emptyset)=0$. 

\begin{proof}
 $P_\D (\emptyset)= P_\D  ( \Omega^C) = \mathfrak p - P(\Omega) =0$. 
\end{proof} 

\item[{\bf (III)}]  If  $A, \,  B  \in \Sigma$ with   $A \subset B$  then  $P_\D (A) $  and  $P_\D (B)$  are  comparable  with  respect  to  the  partial  order $\preceq$, what is more,  
$$P_\D (A) \preceq P_\D (B) . $$

\begin{proof}  
$$ \begin{array}{rcl}
B & = &  B \cap \Omega = B \cap (A \cup A^C) 
\\ & &  \\  & = &  
(B\cap A) \, \cup \, (B \cap A^C) = A \cup (A^C \cap B), 
\end{array}$$ 
where $A \cap (A^C \cap B) = \emptyset$. Hence  $P_\D (B)=P_\D (A) + P_\D (A^C \cap B)$,  and  since  $P_\D (A^C \cap B) \succeq 0$ we can add $P_\D (A) $ to both sides,  proving with this the statement. 
\end{proof}

\end{enumerate}

\begin{Cy}
The  $\D$--probability of any  event is comparable with $\mathfrak p$ and is  $\D$--less or  equal  to  $\mathfrak p$.  
\end{Cy}

Indeed, it is always true that given $A \in \Sigma$, $A \subseteq \Omega$,  hence  $P_\D (A)$ is comparable with  $P_\D (\Omega)$; what is more, $P_\D (A) \preceq P_\D (\Omega) = \mathfrak p$. \\

\begin{Cy}
If  $P_\D(\Omega) = \e$ then for any random event  $A$  there holds that  $P_\D(A) $ is of the form  $\lambda \e$ with $\lambda \in [ 0 , 1]$.    If  $P_\D(\Omega) = \edag$ then for any random event  $A$  there holds that  $P_\D(A) $ is of the form  $\mu \edag$ with $\mu \in [ 0 , 1]$. 
\end{Cy}

\begin{enumerate}
\item[{\bf (IV)}]  The addition theorem.  Given a collection of events $A_1 , \ldots, A_n$, there holds:
$$ \begin{array}{l}
 \displaystyle  P_\D  \left(  \bigcup_{i=1}^n A_n \right) =  \sum_{i=1}^n P_\D (A_i)   -    \sum_{1 \leq i < j \leq n } P_\D (A_i \cap A_j) \;  +
 \\  \\    \quad \quad   + \displaystyle   \sum_{1 \leq i < j < k \leq n } P_\D (A_i \cap A_j \cap A_k) \, +   \;  \cdots + (-1)^{n-1} P_\D (A_1 \cap \cdots \cap A_n) .    
\end{array} $$

\begin{proof}  By induction. For $n=2$
since  $A_1 \cup A_2  =  A_1   \cup  (A_1^C \cap A_2) $ and  also  $A_2 =  (A_1 \cap A_2 )  \cup (A_1^C \cap A_2)$, then  
$$\begin{array}{rcl}
P_\D  (A_1  \cup A_2 ) & = &   P_\D (A_1 ) + P_\D (A_1^C \cap A_2) 
\\ & & \\ & = &  P_\D (A_1) + P_\D (A_2) - P_\D (A_1 \cap A_2).
\end{array}$$  
\end{proof}

\item[{\bf (V)}]  Given two events $A$ and $B$,  $P_\D (A \cup B)$ is comparable  with  $P_\D (A) + P_\D (B)$  and  
$$  P_\D (A \cup B) \preceq  P_\D (A) + P_\D (B).$$
More generally, given  events  $A_1,  \ldots  ,  A_n$ there follows:
$$  \displaystyle  P_\D  \left(  \bigcup_{i=1}^n A_i \right) \preceq   \sum_{i=1}^n P_\D  (A_i). $$

\medskip

\item[{\bf (VI)}]  Theorem of continuity of the $\D$--probability.  

If  $A_1 \supset A_2 \supset \cdots \supset A_n \supset \cdots $  and  $A:=  A_1 \cap A_2 \cap \cdots \cap A_n \cap \cdots  $  then 
$$ \displaystyle  \lim_{n \to \infty} P_\D (A_n) = P_\D (A) = P_\D  \left( \bigcap_{n=1}^\infty A_i \right) .$$  

\begin{proof}  $$ \displaystyle A_n = A \cup  \left( \bigcup_{k=n}^\infty A_k \cap A_{k+1}^C \right) $$
and the summands are  pairwise disjoint, hence
$$ \displaystyle  P_\D (A_n) =P_\D (A) + \sum_{k =n}^\infty P_\D (A_k \cap A_{k+1}^C) .$$
The series here converges for any $n$, in particular,  for $n=1$, hence
$$ \displaystyle  P_\D (A_1) =P_\D (A) + \sum_{k =1}^\infty P_\D (A_k \cap A_{k+1}^C) , $$
thus, the following  sums   go to zero:
$$ \displaystyle \lim_{n \to \infty}   \sum_{k =n}^\infty  P_\D (A_k \cap A_{k+1}^C)  =0  .$$
Finally,  $ \displaystyle  \lim_{ n \to \infty}  P_\D (A_n) = P_\D (A).$ 
\end{proof}

\end{enumerate}

\bigskip

Note that  the  convergence  here is considered with respect  to the hy\-per\-bo\-lic--valued  modulus  $| \cdot |_\k$,  see again \cite{uno}  for  the  details. \\

\begin{Cy} If  $A_1 \subset A_2 \subset \cdots \subset A_n \subset \cdots  $ \, and  $\displaystyle  A:= \bigcup_{n=1}^\infty A_n$, then
$$ \displaystyle   \lim_{n \to \infty}   P_\D (A_n) =  P_\D (A) = P_\D  \left(  \bigcup_{n=1}^\infty A_n \right) .  $$ 
\end{Cy}

\begin{proof}  Take  $B_n:= A_n^C $ and  use  property  ({\bf VI}).
\end{proof}

\section{Conditional probability}

\begin{Dn}
Let  $ ( \Omega , \Sigma , P_\D ) $  be  a  probabilistic space, let $A$ and $B$ be  two events.  The conditional probability  $P_\D (A | B)$ of the  event $A$ under the condition  that $B$ has happened is defined as:

\begin{enumerate}
\item[(1)]  $ \displaystyle   P_\D (A | B) := \frac{ P_\D (A \cap B)}{ P_\D (B)}  $  if  $P_\D (B) \succ 0 $ and  $P_\D (B) \not\in \mathfrak S_\D $;  \\

\item[(2)]  $ \displaystyle  P_\D (A | B) := P_\D (A) $  if  $P_\D (B)=0$;  \\ 

\item[(3)]  $ \displaystyle   P_\D (A | B) := \frac{ P_\D (A \cap B)}{ \lambda_1 } \e  + P_\D (A) \edag    $  if  $P_\D (B)  = \lambda_1 \e $,  $\lambda_1 >0$;    \\

\item[(4)]  $ \displaystyle   P_\D (A | B) := P_\D (A)  \e  +  \frac{ P_\D (A \cap B)}{ \lambda_2 }  \edag   $   if  $P_\D (B) = \lambda_2 \edag$,  $\lambda_2 >0$.   \\

\end{enumerate}
  
\end{Dn}

Let us show that items (3) and (4)  are in a complete agreement with (1). Indeed, using (\ref{cartesian and idempotent representation})   write for any event $A$:  $P_\D (A) = P_1 (A) \e  + P_2 (A) \edag$. Hence,   item (1) in idempotent representation reads:
$$  \displaystyle  P_\D (A | B) = \frac{ P_1 (A \cap B) }{ P_1 (B) }  \e  +   \frac{ P_2 (A \cap B) }{ P_2 (B) }  \edag  =  P_1 (A |B) \e  + P_2 (A |B) \edag  ,$$
meanwhile  items (3) and (4) read:
$$ \begin{array}{rcl}
\displaystyle  P_\D (A | B) &  =  &   \displaystyle    \frac{ P_\D  (A \cap B) }{ \lambda_1 }  \e  +    P_\D  (A)    \edag 
\\  & &  \\  &  = & 
\displaystyle  \frac{   P_1 (A \cap B) \e  + P_2 (A \cap B) \edag }{ P_1 (B) } \e  +  \left(  P_1 (A) \e  + P_2 (A) \edag \right) \edag 
\\  & &  \\  & = & 
\displaystyle  \frac{ P_1( A \cap B) }{ P_1 (B) }  \e  + P_2 (A) \edag  
\\  &  &  \\  & = & 
\displaystyle  P_1 (A |B) \e + P_2 (A |B) \edag 
\end{array}  $$
and  
$$ \begin{array}{rcl}
\displaystyle  P_\D (A | B) &  =  &   \displaystyle    P_\D  (A)    \e  +      \frac{ P_\D  (A \cap B) }{ \lambda_2 }  \edag     
\\  & &  \\  &  = & 
\displaystyle    \left(  P_1 (A) \e  + P_2 (A) \edag \right) \e  +    \frac{   P_1 (A \cap B) \e  + P_2 (A \cap B) \edag }{ \lambda_2 } \edag   
\\  & &  \\  & = & 
\displaystyle   P_1 (A) \e  +     \frac{ P_2( A \cap B) }{ P_2 (B) }  \edag  \; = \;      P_1 (A |B) \e + P_2 (A |B) \edag  .
\end{array}  $$
Thus,  we see a complete compatibility of  the formula in item (1)  and  of    its analogues in items (3) and (4).  \\

Let us show that for a fixed $B$, with $P_\D (B) \not= 0 $,   the conditional probability verifies all the axioms of the $\D$--probability, that is, it defines a new $\D$--probabilistic measure on the measurable  space  $(B,\Sigma_B)$ where $\Sigma_B$ is the $\sigma$--algebra of the sets of the form $A \cap B$ with $A \in \Sigma$.  \\

Indeed, clearly  $P_\D(A|B)  \succeq 0$. Next,  let's see that  $P_\D (B| B) = \mathfrak p$. Indeed:
\begin{enumerate}

\item[(1)]  If  $P_\D (B) \not\in \mathfrak S_\D$, then 
$$  \displaystyle  P_\D (B | B) = \frac{ P_\D (B \cap B)}{ P_\D (B)}   = \frac{P_\D (B)}{ P_\D (B)} =1 . $$

\item[(2)]   If  $P_\D (B) = \lambda_1  \e$,    then 
$$  \displaystyle  P_\D (B | B) = \frac{ P_\D (B \cap B)}{ \lambda_1}  \e  + P_\D (B)   \edag    = \frac{P_\D (B)}{ \lambda_1 }  \e = \e . $$

\item[(3)]   If  $P_\D (B) = \lambda_2  \edag$,    then 
$$  \displaystyle  P_\D (B | B) =  P_\D (B)   \e +   \frac{ P_\D (B \cap B)}{ \lambda_2}  \edag      = \frac{P_\D (B)}{ \lambda_2 }  \edag = \edag . $$

\end{enumerate}
  
Finally, if  $\displaystyle  A = \bigcup_{k=1}^\infty  A_k$  with  $A_i \cap A_j  = \emptyset $ for $i \not= j$ then:  
\begin{enumerate}
\item   If  $P_\D (B) \not\in \mathfrak S_\D$, then 
$$ \begin{array}{rcl}  
\displaystyle  P_\D(A | B)  &  =   &   \displaystyle \frac{P_\D(A \cap B) }{ P_\D(B) }  =  \frac{ P_\D \left(  \bigcup_{k=1}^\infty  A_k  \cap B \right) }{P_\D(B) }
\\  &  &  \\  & = & 
\displaystyle \frac{ \sum_{k=1}^\infty  P_\D (A_k \cap B)}{ P_\D (B)}  = \sum_{k=1}^\infty \frac{ P_\D (A_k \cap B)}{P_\D (B)} 
\\  &  &  \\  & = & 
\displaystyle   \sum_{k=1}^\infty P_\D  (A_k | B)  . 
\end{array}  $$

\item If  $P_\D (B) = \lambda_1 \e \in \mathfrak S_\D$,   since  $A_k \cap B \subset B$ for any  $k$ and since  $A \cap B \subset B$,   write  $P_\D ( A_k \cap B) = \nu_k \e $  and  
$$ \begin{array}{rcl}
P_\D (A \cap B)  &  =  &  \displaystyle    \nu \e = P_1 (A \cap B) \e  =  P_1 \left( \bigcup_{k=1}^\infty  A_k  \cap B \right) \e 
\\  & &  \\  & = &
\displaystyle  \sum_{n=1}^\infty  P_1 (A_k \cap B) \e =   \sum_{n=1}^\infty  \nu_k \e ,
\end{array}  $$ 
hence:
$$ \begin{array}{rcl}
\displaystyle  P_\D(A | B) &  =  &  \displaystyle  \frac{P_\D  (A \cap B) }{ \lambda_1 } \e  + P_\D(A) \edag  =  \frac{ \nu}{\lambda_1 } \e  + P_2(A) \edag 
\\  & &  \\  & = &
\displaystyle  \frac{1}{\lambda_1 }   \sum_{k=1}^\infty \nu_k \e  +   \sum_{k=1}^\infty P_2 (A_k) \edag 
\\  & & \\ &  = & 
\displaystyle   \sum_{k=1}^\infty \left(  \frac{\nu_k}{\lambda_1}  \e + P_2 (A_k) \edag \right)  
\\  & & \\ &  = & 
\displaystyle   \sum_{k=1}^\infty \left(  P_1 (A_k |B) \e + P_2 (A_k ) \edag \right)
\\  & & \\ &  = & 
\displaystyle   \sum_{k=1}^\infty P_\D (A_k |B ) . 
\end{array}$$ 

\item  Similarly if   $P_\D (B) = \lambda_2 \edag \in \mathfrak S_\D$.
\end{enumerate}
Hence, $(B,\Sigma_B, P_\D ( \cdot | B) )$  is  a  new  probabilistic space. \\
\begin{Tm}\label{mult_theorem}(Multiplication Theorem)  
Let $(\Omega, \Sigma, P_\D)$ be a probabilistic  space; let $A$ and $B$  be   two  events. Then
\begin{equation}\label{mult_formula}
P_\D (A \cap B)=P_\D (B) P_\D (A | B) .
\end{equation}
\end{Tm}

\begin{proof}
It is necessary to consider the different cases that arise.

\begin{enumerate}
\item[(a)]  If   $P_\D (B) \succ 0$  and  $P_\D (B) \not\in \mathfrak S_{\D, 0}$, then we know that
$$ \displaystyle  P_\D (A | B) = \frac{ P_\D (A \cap B) }{ P_\D (B) } ,$$
hence  (\ref{mult_formula}) follows. \\

\item[(b)]  If   $P_\D(B) =0$, since  $A \cap B \subset B$  then  $P_\D (A \cap B) =0$   implying  (\ref{mult_formula}). \\

\item[(c)]   If  $P_\D (B) = \lambda_1 \e $ with $\lambda_1 >0$ then
$$ \displaystyle   P_\D (A | B)  =  \frac{  P_\D (A \cap B) }{ \lambda_1 }  \e  + P_\D (A) \edag , $$
hence
\begin{equation}\label{mult2ast}
P_\D(B) P_\D (A | B) = \lambda_1 P_\D (A | B) \e = P_\D ( A \cap B) \e . 
\end{equation}
Since  $A \cap B \subset B$, then  
$$   P_\D (A \cap B)  = P_1 (A \cap B) \e + 0  \edag =P_1 (A \cap B) \e . $$
Rewriting  (\ref{mult2ast})  one gets:
$$ P_\D (B) P_\D (A | B)  = P_1 (A \cap B) \e =  P_\D (A \cap B).   $$

\item[(d)]     If  $P_\D (B) = \lambda_2 \edag  $ with $\lambda_2 >0$, one proceeds as in  (c).  
\end{enumerate}
\end{proof}

This theorem has a generalization for $n$  random events.

\begin{Tm} (Generalized multiplication theorem.)  Let  $A_1 ,  \ldots , A_n$ be random  events.  If  they  satisfy  any  of the  following  conditions:

\begin{enumerate}

\item[(1)]     $P_\D (A_1 \cap \cdots \cap A_n)$  is not  a  zero--divisor.  \\

\item[(2)]   \begin{enumerate}
\item[(a)]  There exists $k_0   \in  \{ 1 , \ldots , n \}  $  such that  $P_\D (A_{k_0})   = \mu_{k_0} \e$, with  $\mu_{k_0} >0$, i.e.,   $ P_\D (A_{k_0} ) $ is a   zero--divisor  in   $\D_\e^+$,    \\
\end{enumerate}
\noindent and
\begin{enumerate}
\item[(b)]  $ \displaystyle  P_\D \left( \bigcap_{\ell =1}^n A_\ell \right) $  belongs  to  $\D_\e^+$  also.  \\

\end{enumerate}  

\item[(3)]   \begin{enumerate}
\item[(a)]  There exists $k_0   \in  \{ 1 , \ldots , n \}  $  such that  $P_\D (A_{k_0})   = \lambda_{k_0} \edag$, with  $\lambda_{k_0} >0$, i.e.,   $ P_\D (A_{k_0} ) $ is a   zero--divisor  in   $\D_{\edag}^+$,    \\
\end{enumerate}
\noindent and
\begin{enumerate}
\item[(b)]  $ \displaystyle  P_\D \left( \bigcap_{\ell =1}^n A_\ell \right) $  belongs  to  $\D_{\edag}^+$  also;

\end{enumerate}  

\end{enumerate}
then  
\begin{equation}\label{form_theor_mult}  
\displaystyle  P_\D  \left(  A_1 \cap \cdots \cap A_n   \right)   = P_\D(A_1) P_\D(A_2 | A_1) \cdots P_\D(A_n | A_1 \cap \cdots \cap A_{n-1} ) . 
\end{equation}

\end{Tm}

\begin{proof}  Assume  that  (1)  occurs.  Since
\begin{equation}\label{contenciones}  
\displaystyle  \bigcap_{i=1}^{n-1} A_i \subset \bigcap_{i=1}^{n-2} A_i \subset \cdots \subset A_1
\end{equation}
the hypothesis $\displaystyle  P_\D \left( \bigcap_{\ell=1}^n A_\ell \right)   =  P_1 \left( \bigcap_{\ell=1}^n A_\ell \right)  \e  +   P_2 \left( \bigcap_{\ell=1}^n A_\ell \right)  \edag \not\in  \mathfrak S_{\D, 0}$  is  equi\-va\-lent to say  that    $\displaystyle P_1 \left( \bigcap_{\ell=1}^n A_\ell \right)  >0$  and   $\displaystyle P_2 \left( \bigcap_{\ell=1}^n A_\ell \right)  >0$, and  it  implies  that      $P_\D (A_\ell) \not\in \mathfrak S_{\D, 0}$  for  all   $\ell \in \left\{ 1, \ldots , n \right\} $  and  also  that  
$\displaystyle  P_\D \left( \bigcap_{i=1}^{n-k} A_i \right) $ is a  strictly  positive   hyperbolic  number  (that is positive  and  not  zero--divisor)  for all  $k \in  \{1, \ldots , n-1 \} $ and thus all conditional  probabilities of the form  $\displaystyle  P_\D \left( A_k \big| \bigcap_{i=1}^{k-1} A_i \right) $ for  $ k \in \{2, \ldots , n \} $ are well--defined, implying that
$$ \begin{array}{rcl}
\displaystyle    P_\D \left( \bigcap_{\ell=1}^n A_\ell \right)   &  =  &  \displaystyle    P_1 \left( A_1  \right) P_1 \left(  A_2 \mid A_1  \right) \cdots P_1 \left(  A_n \mid  \bigcap_{\ell=1}^{n-1}  A_\ell \right)  \,   \e      \; +
\\  &  &   \\  &  &  
\displaystyle   \qquad  + \;     P_2 \left( A_1  \right) P_2 \left(  A_2 \mid A_1  \right) \cdots P_2 \left(  A_n \mid  \bigcap_{\ell=1}^{n-1}  A_\ell \right)  \,   \edag , 
\end{array}  $$
hence (\ref{form_theor_mult}) follows.   

Assuming  now  that  (2)  occurs, the hypothesis  that   $  \displaystyle    P_\D \left( \bigcap_{\ell=1}^n A_\ell \right)  $  is  a  po\-si\-ti\-ve  zero--divisor,  let's say  $ \displaystyle    P_\D \left( \bigcap_{\ell=1}^n A_\ell \right)  = \mu \e $,  $\mu >0$,  implies    that   $ \displaystyle    P_1 \left( \bigcap_{\ell=1}^n A_\ell \right)   = \mu >0$.  

On the  other  hand,  assume  that  $k_0$ is the minimum integer in  $ \left\{ 1, 2, \ldots , n \right\}$  such  that   $  \displaystyle  P_\D  \left( A_{k_0} \right) $  is  a  zero--divisor.  This  implies  that $  \displaystyle    P_\D \left( \bigcap_{\ell=1}^k A_\ell \right)  \in  \D_\e $  for  all  $k \geq k_0$;  writing  in this case  $\displaystyle    P_\D \left( \bigcap_{\ell=1}^k A_\ell \right)  = \nu_k \e$, one has:
$$  \displaystyle     P_\D \left(   A_{k+1} \big|   \bigcap_{\ell=1}^k A_\ell \right)  =  \frac{    P_\D \left( \displaystyle  \bigcap_{\ell=1}^{k+1}  A_\ell \right) }{ \nu_k } \,   \e  +  P_\D \left(  A_{k+1} \right) \, \edag \qquad {\rm for} \quad k \geq k_0 . $$  
We  have  also that  
$$ \begin{array}{rcl}
 \displaystyle     P_\D \left(   A_{k_0} \big|   \bigcap_{\ell=1}^{k_0 -1}   A_\ell \right)  &  =  &  \displaystyle     \frac{    P_\D \left( \displaystyle  \bigcap_{\ell=1}^{k_0}  A_\ell \right) }{  P_\D \left( \displaystyle  \bigcap_{\ell=1}^{k_0-1}  A_\ell \right) }  \, =
 \\   &  &   \\  &  = &
 \displaystyle      \frac{\nu_k \e }{    P_1 \left( \displaystyle  \bigcap_{\ell=1}^{k_0-1 }  A_\ell \right)   \, \e  +   P_2 \left( \displaystyle  \bigcap_{\ell=1}^{k_0-1 }  A_\ell \right)   \, \edag  }  
 \\  &  &  \\   &  =  &    \displaystyle       \frac{\nu_k  }{    P_1 \left( \displaystyle  \bigcap_{\ell=1}^{k_0-1 }  A_\ell \right) }    \, \e  .  
 \end{array}   $$
 Hence
$$ \begin{array}{l}
\displaystyle   P_\D \left( A_1  \right)  P_\D \left( A_2 \big| A_1  \right)  \cdots    P_\D \left( A_{k_0}  \big|    \bigcap_{\ell=1}^{k_0-1 }  A_\ell   \right)   \cdots    P_\D \left( A_n  \big|    \bigcap_{\ell=1}^{n-1 }  A_\ell   \right)   \,  = 
\\  \\ 
\displaystyle  \; = \;  P_1 \left( A_1  \right)  P_1 \left( A_2 \big| A_1  \right)  \cdots    P_1 \left( A_{k_0}  \big|    \bigcap_{\ell=1}^{k_0-1 }  A_\ell   \right)   \cdots    P_1 \left( A_n  \big|    \bigcap_{\ell=1}^{n-1 }  A_\ell   \right)   \, \e    
\\  \\ 
\displaystyle  \; = \;  P_1 \left( A_1 \cap \cdots \cap A_n  \right)  \, \e = P_\D \left( A_1  \cap \cdots \cap A_n  \right) .
\end{array}  $$

\medskip

The case (3)  is  proved  analogously.  
\end{proof}

\medskip

\begin{Dn}
Let    $A$ and $B$  be  two   random  events.

\begin{enumerate}

\item  $A$ is called  independent  of   $B$  if  
$$  P_\D (A | B) = P_\D (A). $$

\item  $B$  is  called  independent  of   $ A$  if  
$$  P_\D (B | A) = P_\D (B). $$

\item  $A$  and  $B$  are  called  mutually  independent     if    $A$  is independent of  $B$  and  $B$  is  independent  of  $A$.  \\
    
\end{enumerate}    

\end{Dn}

Let us  analyze  all possible  situations.  

\begin{enumerate}

\item[(i)]  Assume  that  $P_\D(A) = P_\D (B)=0$.  By definition, in this case  $P_\D (A | B) = P_\D (A) $ and  $P_\D(B |A) = P_\D (B)$, thus  $A$ and $B$  are mutually independent. Moreover, it is enough to assume  that  one of  the two  probabilities  only, say  $P_\D (A)$,  equals  zero. Indeed,  if  $P_\D (A)=0$ then  $  P_\D (A \cap B) =0$ and hence   
$$  P_\D (B | A) = P_\D (B) \qquad {\rm and} \qquad  P_\D (A |B) =0=P_\D (A),$$
i.e.,   $A$  and $B$ are  mutually  independent. And in any case there follows $ P_\D (A \cap B) = P_\D (A) P_\D (B)$. \\

\item[(ii)]  Assume  that  both  $P_\D(A)$ and $P_\D(B)$ are not  in $\mathfrak S_{\D , 0}$.     That  $A$ is  independent  of  $B$ is equivalent to  
$$ \displaystyle  P_\D(A|B) = P_\D(A) = \frac{ P_\D(A \cap B) }{ P_\D (B) } $$
which, in turn, is equivalent to 
$$ P_\D  (A \cap B) =P_\D (A) \cdot P_\D (B).$$
In the same way, if $B$  is independent   of  $A$  then this is equivalent to  
$$P_\D(B |A)  = P_\D(B) = \frac{ P_\D (A \cap B) }{ P_\D(A) } , $$
i.e., 
$$P_\D(A \cap B) = P_\D(B) \cdot P_\D(A).$$
Finally, under the assumed hypotheses $A$ is independent of $B$ if and only if $B$ is independent of $A$  if and only if   $A$  and $B$  are  mutually  independent. \\

\item[(iii)]   Assume that $P_\D(A) $  and  $ P_\D (B)$ are zero--divisors which   both belong to  $\D_\e$:  $ P_\D(A) = \lambda \e $  and $ P_\D (B) = \mu \e$  with positive reals  $\lambda$ and $\mu$.  This means  that  $P_1 (A)= \lambda$ and $P_1 (B) = \mu$  meanwhile  $P_2(A) =P_2(B) =0$  implying  that  $  0  \leq P_1 (A \cap B) =: \nu$  and  $P_2 (A \cap B)=0$  (hence $P_\D (A \cap B) = \nu \e$).  Suppose  that  $A$ is independent of $B$; this is equivalent to  
\begin{equation}\label{probAyBzero_div_e}
\begin{array}{rcl}
\lambda \e & = &  \displaystyle  P_\D(A) = P_\D (A | B)  = \frac{ P_\D (A \cap B)}{ \mu } \e  + P_\D (A) \edag
\\  & & \\  & = &
\displaystyle  \frac{ P_1 (A \cap B) }{ \mu } \e = \frac{ \nu }{ \mu} \e ;
\end{array} 
\end{equation}
thus $A$ is independent of $B$ if and only if  
$$ \displaystyle  \lambda = \frac{ \nu }{ \mu} . $$
Assume the equality  $ \displaystyle  \lambda = \frac{ \nu}{ \mu} $ holds.  Considering now $P_\D (B | A)$  we have:
\begin{equation}\label{probAyBzero_div_e-2}
\begin{array}{rcl}
P_\D (B|A) & = &  \displaystyle  \frac{  P_\D (A \cap B) }{ \lambda } \e  + P_\D (B) \edag = \frac{ \nu }{ \lambda } \e = \frac{ \nu}{ \nu/ \mu} \e
\\  & &  \\  & = &
\displaystyle   \mu \e =P_\D (B) 
\end{array} 
\end{equation}
which  means that $B$ is independent  of $A$,  and  thus $A$  and $B$  are mutually independent.  Observe that using  (\ref{probAyBzero_div_e})  or  (\ref{probAyBzero_div_e-2}) one concludes that  for  independent  $A$  and  $B$  
$$ P_\D (A \cap B) = P_\D (A \cap B) \e =  \lambda  \mu \e = (\lambda \e)  (\mu \e) = P_\D (A) P_\D(B) .$$ 

In the same way the case of $P_\D(A) $ and $ P_\D (B)$ being  both in $\D_\edag$ is covered. \\

\item[(iv)]  Assume that  $P_\D(A) $ and $P_\D(B)$  are zero--divisors  but  now such that  $P_\D(A) = \lambda \e \not= 0$  and $P_\D (B) = \mu \edag \not=0,$  or vice--versa.  This gives, in particular,  that    since  $A \cap B \subset A$  and  $A \cap B \subset B$  then  
$$  P_1 (A \cap B) =P_2 (A \cap B)=0  ,  $$
that is   $P_\D (A \cap B) =0$.  
Consider  $P_\D (A | B)$,   one has:
$$ \displaystyle  P_\D (A | B) = \frac{  P_\D (A \cap B) }{ \mu } \edag + P_\D (A) \e = \lambda \e =  P_\D (A) ,$$
which means  that the hypotheses imply that  $A$  is independent of $B$.  But   
$$ \displaystyle  P_\D ( B | A) = \frac{ P_\D (A \cap B) }{ \lambda } \e  + P_\D (B) \edag = P_\D (B) \edag =  \mu \edag =  P_\D (B),$$
that is, $B$  is independent of $A$. 

Somewhat  paradoxically,  in this case  $A$ and $B$  are always mutually  independent. And there follows that
$$ P_\D (A) P_\D (B) = ( \lambda \e )(\mu \edag ) =0= P_\D (A \cap B). $$

\item[(v)]  The last case assumes that $P_\D(A) $ is  zero--divisor  and $P_\D (B) $ is an invertible  hyperbolic   number  or  vice--versa.  Set  
$$  P_\D (A)  = \lambda \e \not= 0,  \quad    P_\D (B) =  \mu_1 \e  + \mu_2 \edag \not\in \mathfrak S_{ \D , 0} . $$   
In particular, this implies that  
$$ \nu:=  P_1  ( A \cap B) \geq 0 \quad  {\rm and }  \quad  P_2 (A \cap B)=0   .  $$  
That $A$ is  independent of $B$  is equivalent  to  
$$ \begin{array}{rcl}
P_\D(A |B) & = &  \displaystyle   \frac{P_\D ( A \cap B) }{ P_\D (B) } = \frac{  \nu \e }{ \mu_1 \e  + \mu_2 \edag } =  \frac{ \nu }{ \mu_1 } \e
\\  & &  \\  & = & 
\displaystyle  P_\D(A) = \lambda \e ,
\end{array} $$ 
i.e.,  $ \displaystyle  \lambda = \frac{ \nu }{ \mu_1} $.  

On the other hand,
$$ \begin{array}{rcl}
P_\D (B |A) & =  &  \displaystyle   \frac{ P_\D ( A \cap B) }{ \lambda }  \e  +  P_\D (B) \edag  =  \frac{P_1(A \cap B)}{\lambda} \e + P_2 (B) \edag 
\\ & & \\ & = &
\displaystyle  \frac{ \nu }{ \lambda } \e + P_2 (B) \edag  = \mu_1 \e + \mu_2 \edag  = P_\D (B) ,
\end{array} $$
i.e.,  $B$ is independent  of $A$. Of course, the reasoning is reversible, thus $A$ and $B$ are mutually  independent  if and only if  one of them is independent of the another one.  

Finally note that from the above one has for independent events  that  $  \nu =  \lambda \mu_1$, hence:
$$  P_\D (A \cap B)= \nu \e = (\lambda \e )(\mu_1 \e) = P_\D (A) P_\D(B) .$$
  
\end{enumerate}

\medskip

Resuming the  just made  analysis , one has the following

\begin{cor} Given two random events $A$ and $B$, then $A$ is independent of $B$ if and only if $B$ is independent of $A$.  \\
\end{cor}

\begin{cor} If $A$ and $B$ are  mutually  independent events then the multiplication  theorem  becomes
\begin{equation}\label{independent-multiplication}
P_\D (A \cap B) = P_\D (A) P_\D(B) .
\end{equation}
\end{cor}

\medskip

\begin{Tm}
If $A$ and $B$ are  mutually  independent events then so are $A$ and $B^C$,  $A^C$  and  $B$,  $A^C$ and   $B^C$.
\end{Tm}

\begin{proof}  It is enough to prove  for $A$ and  $B^C$.  Write  $A =  (A \cap B) \cup  (A \cap B^C)$ with  disjoint summands. Then  $P_\D(A) = P_\D(A \cap B) + P_\D (A \cap B^C)$  from  where  
\begin{equation}\label{AyBcompl}
P_\D (A \cap B^C)   =       P_\D(A) -  P_\D(A \cap B) =  P_\D(A)  -  P_\D(A) P_\D(B). 
\end{equation}
Of  course,  we  can  factorize  $P_\D(A)$  but  the consequent computation    depends  on  the  value  of  $P_\D(\Omega)$.  Thus    we have to consider the following cases:

\begin{enumerate}

\item[(1)]  $P_\D(A) \not\in \mathfrak S_{\D, 0}$, then  $P_\D(\Omega) \not\in \mathfrak S_{\D, 0}$ and in this case $\mathfrak p =1$. Hence
$$ P_\D (A\cap B)= P_\D (A) (1 - P_\D (B) ) = P_\D (A)  P_\D (B) ,$$
thus $A$ and $B^C$ are mutually independent.

\item[(2)] If $P_\D (A) \in \mathfrak S_{\D, 0}$, let us say,  $P_\D (A)= \lambda \e$, since $A \cap B^C \subset A$, then  $P_\D (A \cap B^C) = \nu \e$ for some $\nu \geq 0$; write $P_\D (B)=\mu_1 \e + \mu_2 \edag$,   and one has two more subcases:

\begin{enumerate}

\item[(a)] If $P_\D (\Omega) =1 \not\in \mathfrak S_{\D, 0}$, then  $P_\D(B^C)= 1 - P_\D(B) = (1 - \mu_1 )\e + (1- \mu_2) \edag $ and there follows:
$$ \begin{array}{rcl}
\nu \e & = &  P_\D(A \cap B^C)
\\ & & \\  & = &  
\displaystyle   \lambda \e - \lambda \e \left(  \mu_1 \e + \mu_2 \edag \right) \, = \,  \lambda \e - \lambda \e (  \mu_1 \e )
\\  & &  \\  & = & \displaystyle
\lambda \left( 1 - \mu_1 \right) \e =  \left( \lambda \e \right) \left(  (1 - \mu_1 ) \e \right) = P_\D (A) P_\D \left( B^C \right) , 
\end{array} $$
thus, $A$ and $B^C$ are mutually independent.  \\

\item[(b)]  If $P_\D (\Omega)= \e$ (note that the  equality  $P_\D(\Omega)= \edag$ is impossible)   then  ne\-ce\-ssa\-ri\-ly  $P_\D (B) = \mu_1 \e$ and   $P_\D (B^C) = (1 - \mu_1)\e$, hence
$$ \begin{array}{rcl}
P_\D (A \cap B^C) &  =  &    \lambda \e - \lambda \e \mu_1 \e
\\ & &  \\  &  = &  (\lambda \e)( 1 - \mu) \e =  P_\D (A) P_\D (B^C) ,
\end{array} $$
and $A$ and $B^C$ are mutually independent. \\

\end{enumerate}

\item[(3)] The case $P_\D(A) = \lambda_2 \edag $ is treated similarly.

\end{enumerate}

\end{proof}

\begin{Dn}
Random  events  $A_1 , \ldots , A_n$  are called  mutually  (or jointly) independent if  for  any  subset  of  indices  $i_1 , \ldots , i_r$  with  $ 1 \leq i_1 < i_2 <  \cdots < i_r \leq n$ ($ r \in \{ 2, \ldots n \}   $)   there  holds:
$$  P_\D(A_{i_1} \cap  \cdots \cap  A_{i_r} )= P_\D(A_{i_1} ) \cdots  P_\D (A_{i_r}) . $$
If this holds  for $r=2$ only then the events  are  called  pair--wise  independent.  Generally  speaking,  pair--wise  independence  and  joint  independence  are  different  notions.  \\
\end{Dn}

If   $A_1 , \ldots , A_n$  are  mutually independent  events  then the general  multiplication  theorem  holds in a  simplified form:
$$  P_\D(A_1 \cap \cdots \cap A_n)=P_\D(A_1)  \cdots  P_\D(A_n).  $$

\medskip

\begin{Dn}
Let $(\Omega, \Sigma ,P_\D)$ be  a  $\D$--probabilistic  space. Let  $H_1, \ldots , H_n$  be  pairwise  disjoint random  events  with  (not necessarily  strictly)  positive  probabilities and such  that  $H_1 \cup \cdots \cup H_n =\Omega$.   Then  the collection  $\{ H_1 , \ldots , H_n \}$ is called  a  fundamental  (or  complete)  system of  events  (FSE).  \\
\end{Dn}

\begin{Tm}
(hyperbolic law of total probability; complete hyperbolic  probability  formula). Given  $(\Omega , \Sigma , P_\D)$, $A$  a  random  event;  $\{ H_1 , \ldots , H_n \}$    a  FSE, then 
$$ \displaystyle   P_\D (A) = \sum_{i=1}^n  P_\D (H_i) P_\D(A | H_i) ;$$
\end{Tm}

\begin{proof}  Since  $ \displaystyle  A=A \cap \Omega= A \cap \left( \bigcup_{i=1}^n  H_i \right)  =  \bigcup_{i=1}^n  A \cap H_i $ and the  events  $ A \cap H_i $  are  pairwise  disjoint then
$$ \displaystyle   P_\D(A) =\sum_{i=1}^n  P_\D \left( A \cap H_i \right) .  $$
It is enough to apply  now    Theorem \ref{mult_theorem} (``Multiplication Theorem"). 
\end{proof}

\begin{Tm} (Bayes' theorem).  Let $(\Omega, \Sigma, P_\D)$, $A$  and  $\{ H_1 , \ldots , H_n \}$ be  as in   the previous  theorem, then:

\begin{enumerate}

\item[1)]  if      $P_\D(A)$  is  an  invertible  hyperbolic  number  then  

\begin{equation}\label{estrella}  
\displaystyle  P_\D (H_k|A) =  \frac{  P_\D (H_k) \cdot P_\D (A | H_k)}{ \sum_{i=1}^n  P_\D (H_i) \cdot  P_\D (A | H_i) } =  \frac{  P_\D (H_k) \cdot P_\D (A | H_k)}{ P_\D (A) } ; 
\end{equation}

\item[2)]  if $P_\D(A)  = \lambda \e $  with $\lambda >0$ then     

\begin{equation}\label{dos-estrellas}  
\displaystyle \left( P_\D (H_k) \cdot  P_\D (A | H_k ) -   P_\D (H_k | A) \cdot  \sum_{i=1}^n  P_\D (H_i) \cdot  P_\D (A | H_i) \right) \e =0.
\end{equation}

\item[3)]  If  $P_\D(A)  = \mu \edag $  with $\mu >0$ then     

\begin{equation}\label{tres-estrellas}  
\displaystyle \left( P_\D (H_k) \cdot  P_\D (A | H_k ) -   P_\D (H_k | A) \cdot  \sum_{i=1}^n  P_\D (H_i) \cdot  P_\D (A | H_i) \right) \edag =0.
\end{equation}

\end{enumerate}
\end{Tm}

\begin{proof}   

\begin{enumerate}

\item[1)]  Let  $P_\D (A)$  be  an invertible hyperbolic number then  the multiplication theorem  gives:
$$P_\D (A \cap H_k) = P_\D (H_k)  P_\D (A | H_k) = P_\D (A) P_\D  (H_k |A) .  $$
Thus, a part  of formula  (\ref{estrella})  verifies. Using the  hyperbolic  law of total probability   gives the  rest  of   (\ref{estrella}). \\

\item[2)]  Let  $P_\D (A) = \lambda \e $   with $\lambda $  being a positive real number then the multiplication theorem  leads to the equality: 
\begin{equation}\label{ultima}
P_\D(A)  P_\D (H_k | A)  = P_\D (H_k)  P_\D (A | H_k) . 
\end{equation}
Note that the left--hand side of  (\ref{ultima})  is an element of  $\D_\e$, hence the right--hand side  must be an element of $\D_\e$ also. But this is true indeed, since  in the definition of  $P_\D (A | H_k)$  the factor  $P_\D (A \cap H_k)$ is involved and because  $A \cap H_k \subset A$, hence $P_\D (A \cap H_k) \in \D_\e$. Now, recalling the property of Remark \ref{useful property}, equation (\ref{ultima})  can be rewritten as
$$ \begin{array}{rcl}
0 & = &  P_\D (H_k)  P_\D (A | H_k)  -   P_\D (H_k | A) P_\D(A)
\\  & &  \\ & = & 
\displaystyle \left( P_\D (H_k)  P_\D (A | H_k)  -   P_\D (H_k | A) P_\D(A) \right) \e ,
\end{array} $$
finally using the hyperbolic  law of total probability one obtains (\ref{dos-estrellas}). \\

\item[3)]  We proceed similarly to item 2).
\end{enumerate}
\end{proof}

\bigskip


\bigskip

\bigskip

\end{document}